\documentclass[11pt]{article}
\usepackage[numbers,square]{natbib}
\usepackage{amsmath}
\usepackage{amsfonts}
\usepackage{amssymb}
\usepackage{commath}
\usepackage{graphicx}
\usepackage{xspace}
\usepackage{color}
\usepackage[lofdepth,lotdepth]{subfig}
\usepackage{stmaryrd}
\usepackage{hyperref}
\usepackage{cleveref}
\usepackage{url}
\usepackage{amsthm}
\usepackage{footnote}
\makesavenoteenv{tabular}
\makesavenoteenv{table}
\usepackage{tikz}

\usepackage[margin=0.7in]{geometry}

\newcommand{\jump}[1]{[\![#1]\!]}
\begin{document}
\begin{table}[tbp]
\begin{center}
\begin{tabular}{c}
  \textsc{A locally conservative and energy-stable finite element
  method}\\\textsc{for the Navier--Stokes problem on time-dependent domains}\\
  \\
  \textsc{Tam\'as L. Horv\'ath\footnotemark\footnotetext{Department of Applied Mathematics, University of Waterloo, Waterloo, Ontario, Canada. Now at Department of Mathematics and Statistics, Oakland University, Rochester, Michigan, USA. Email: \url{thorvath@oakland.edu}} and Sander Rhebergen\footnote{Department of Applied Mathematics, University of Waterloo, Waterloo, Ontario, Canada. Email: \url{srheberg@uwaterloo.ca}. ORCID: 0000-0001-6036-0356.}}
\end{tabular}
\end{center}
\end{table}
\subsubsection*{Abstract}
We present a finite element method for the incompressible
Navier--Stokes problem that is locally conservative, energy-stable and
pressure-robust on time-dependent domains. To achieve this, the
space--time formulation of the Navier--Stokes problem is
considered. The space--time domain is partitioned into space--time
slabs which in turn are partitioned into space--time simplices. A
combined discontinuous Galerkin method across space--time slabs, and
space--time hybridized discontinuous Galerkin method within a
space--time slab, results in an approximate velocity field that is
$H({\rm div})$-conforming and exactly divergence-free, even on
time-dependent domains. Numerical examples demonstrate the convergence
properties and performance of the method.
\\
\\
\textbf{Keywords:} Navier--Stokes, Space--Time, Hybridized,
Discontinuous Galerkin, Time-Dependent Domains.
%
\section{Introduction}
\label{s:introduction}

Space--time discontinuous Galerkin (DG) finite element methods have
proven to be excellent discretization methods for the solution of
partial differential equations on time-dependent domains. For example,
space--time DG methods have successfully been applied to the
compressible Euler~\cite{Vegt:2002} and Navier--Stokes
equations~\cite{Klaij:2006, Wang:2015}, incompressible
flows~\cite{Vegt:2008, Rhebergen:2013b, Tavelli:2015, Tavelli:2016},
shallow water equations~\cite{Ambati:2007}, nonlinear water
waves~\cite{Vegt:2007} and two-phase flows~\cite{Rhebergen:2009,
  Sollie:2011}. The success of space--time discontinuous Galerkin
methods lie in that they automatically satisfy the geometric
conservation law (uniform flow on a dynamic mesh remains
uniform)~\cite{Lesoinne:1996}, they can be made unconditionally
stable, and are fully conservative and higher-order accurate in space
and time. Space--time methods are furthermore well suited for
$hp$-adaptivity in both space and time~\cite{Vegt:2002}.

Space--time DG methods, however, are computationally expensive. A
partial differential equation on a $d$-dimensional time-dependent
domain is discretized by a space--time method in $d+1$ space--time,
adding an extra dimension to the problem. This results in large
systems of (non)linear algebraic equations that need to be
solved. Although solvers exist for these systems of
equations~\cite{Neumueller:2013, Vegt:2012a, Vegt:2012b}, the size of
the problem remains an issue.

The hybridizable discontinuous Galerkin (HDG) finite element method
was introduced in~\cite{Cockburn:2009} to reduce the computational
cost of DG methods. To achieve this, HDG methods are constructed such
that the only globally coupled degrees-of-freedom lie on cell
boundaries. This results in a significant reduction in the size of the
problem compared to standard DG methods where the globally coupled
degrees-of-freedom lie on cell interiors. It was for this reason
that~\cite{Rhebergen:2012, Rhebergen:2013} introduced space--time HDG
methods as computationally cheaper alternatives to space--time DG
methods.

The first space--time HDG method for the incompressible Navier--Stokes
equations was introduced in~\cite{Rhebergen:2012}. Using space--time
hexahedral cells, they were able to obtain optimal rates of
convergence for the velocity, velocity gradient and pressure fields on
time-dependent domains. A drawback of their method, as is common with
many other discontinuous Galerkin methods for incompressible flows, is
that their method cannot both be locally momentum conserving and
energy-stable due to the discrete velocity field not being point-wise
divergence-free and
$H(\text{div})$-conforming~\cite{Cockburn:2005}. The absence of a
velocity field that is point-wise divergence-free and
$H(\text{div})$-conforming also has other consequences, for example,
the velocity error will depend on the pressure error scaled by the
inverse of the viscosity~\cite{John:2017}. This lack of
`pressure-robustness' may cause large errors in the velocity,
especially for convection dominated flows.

Different techniques have been developed for incompressible flows to
obtain DG methods on fixed domains that result in point-wise
divergence-free and $H(\text{div})$-conforming velocity fields, such
as post-processing~\cite{Cesmelioglu:2016} and the use of
$H(\text{div})$-conforming finite elements~\cite{Cockburn:2007,
  Lehrenfeld:2016}. An alternative is the use of HDG methods. In a
series of papers, Rhebergen and Wells introduced HDG finite element
methods for incompressible flows on fixed domains resulting in
discrete velocity fields that are point-wise divergence-free and
$H(\text{div})$-conforming. In~\cite{Rhebergen:2017} they introduced
and analyzed their HDG method for the Stokes problem proving optimal
error estimates. In~\cite{Rhebergen:2018b} optimal preconditioners
were introduced for the fast solution of the resulting linear
systems. Their HDG method was extended in~\cite{Rhebergen:2018a} to
the Navier--Stokes problem resulting in a scheme that is momentum
conserving, energy stable and pressure-robust.

In this paper we present a space--time HDG method for the
Navier--Stokes problem on \emph{time-dependent} domains that results in
point-wise divergence-free and $H(\text{div})$-conforming velocity
fields. We will show that our scheme is momentum conserving, energy
stable and pressure-robust. To the best of our knowledge, this is the
first finite element method that achieves all these properties on
time-dependent domains. To achieve this we divide the whole
space--time domain into space--time slabs. Each space--time slab is
then divided into space--time tetrahedra. This is different from many
other space--time methods in which space--time slabs are typically
divided into space--time prisms or hexahedra. Given the space--time
tetrahedra, the HDG method of~\cite{Rhebergen:2018b} can naturally be
extended to a space--time formulation on time-dependent domains.

The rest of this paper is organized as
follows. In~\cref{s:navierstokes} we introduce the Navier--Stokes
problem after which we introduce the space--time HDG method
in~\cref{s:discretenavierstokes}. We discuss properties of the
space--time HDG method in~\cref{s:properties} and show numerical
results in~\cref{s:numericalexamples}. Conclusions are drawn
in~\cref{s:conclusions}.

\section{The Navier--Stokes problem on time-dependent domains}
\label{s:navierstokes}
Let $\Omega(t) \subset \mathbb{R}^d$ be a time-dependent polygonal
($d=2$) or polyhedral ($d=3$) domain and $I=(0,T)$ the time interval
of interest. We consider the incompressible Navier--Stokes equations
on the space--time domain
$\mathcal{E} := \cbr{ \boldsymbol{x}\in \Omega(t),\ t\in I}$:
\begin{subequations}
  \label{eq:navsto}
  \begin{align}
    \label{eq:navsto_mom}
    \partial_t \boldsymbol{u} + \nabla \cdot \del{\boldsymbol{u} \otimes \boldsymbol{u}} + \nabla p 
    - \nu \nabla^2 \boldsymbol{u} &= \boldsymbol{f}
    && \text{in} \ \mathcal{E},
    \\
    \label{eq:stnavsto_mass}
    \nabla \cdot \boldsymbol{u} &= 0 && \text{in} \ \mathcal{E},
  \end{align}
\end{subequations}
where $\boldsymbol{u}:\mathcal{E}\to\mathbb{R}^d$ is the velocity field,
$p:\mathcal{E}\to\mathbb{R}$ is the kinematic pressure,
$\boldsymbol{f}:\mathcal{E}\to\mathbb{R}^d$ a given forcing term and $\nu$
the constant kinematic viscosity.

The boundary of the space--time domain $\mathcal{E}$ is partitioned
such that
$\partial\mathcal{E} = \partial\mathcal{E}^D
\cup \partial\mathcal{E}^N \cup \Omega(0) \cup \Omega(T)$,
where there is no overlap between any two of the four sets. Here
$\partial\mathcal{E}^D$ and $\partial\mathcal{E}^N$ denote,
respectively, the part of the space--time boundary with Dirichlet and
Neumann boundary conditions. The space--time outward unit normal
vector to $\partial\mathcal{E}$ is denoted by $(n_t, \boldsymbol{n})$, with
$n_t\in\mathbb{R}$ the temporal and $\boldsymbol{n}\in\mathbb{R}^d$ the spatial
component. We then impose the following initial and boundary
conditions:
\begin{subequations}
  \label{eq:boundaryconditions}
  \begin{align}
    \label{eq:st_bc_dirichlet}
    \boldsymbol{u} &= \boldsymbol{0} && \text{on} \ \partial\mathcal{E}^D,
    \\
    \label{eq:st_bc_neumann_outflow}
    \sbr{n_t + \boldsymbol{u}\cdot\boldsymbol{n} - \max\del{n_t + \boldsymbol{u}\cdot\boldsymbol{n}, 0}}\boldsymbol{u} + 
    \del{p\mathbb{I} - \nu\nabla \boldsymbol{u}}\cdot\boldsymbol{n} &= \boldsymbol{g} && \text{on} \ \partial\mathcal{E}^N,
    \\
    \label{eq:st_ic}
    \boldsymbol{u}(0, \boldsymbol{x}) &= \boldsymbol{u}_0(\boldsymbol{x}) && \text{in}\ \Omega(0),
  \end{align}
\end{subequations}
where $\boldsymbol{g}:\partial\mathcal{E}^N \to \mathbb{R}^d$ is given
Neumann boundary data, $\boldsymbol{u}_0:\Omega(0) \to \mathbb{R}^d$ is a
given divergence-free velocity field, and $\mathbb{I}$ is the
$d\times d$ identity matrix.

\section{The discrete Navier--Stokes problem}
\label{s:discretenavierstokes}
In this section we introduce the space--time hybridized discontinuous
Galerkin method for the Navier--Stokes problem~\cref{eq:navsto}
and~\cref{eq:boundaryconditions}.

\subsection{The space--time mesh and the finite element function spaces}
To introduce the space--time mesh, we first partition the time
interval $(0, T)$ using time levels
$0 = t^0 < t^1 < \cdots < t^N = T$. The $n$-th time interval is
defined as $I^n = (t^n, t^{n+1})$, which has length
$\Delta t^n = t^{n+1} - t^n$. Space--time slabs are then defined as
$\mathcal{E}^n := \cbr{(t,\boldsymbol{x})\in\mathcal{E} : t\in I^n}$, which
has boundaries $\Omega^n:=\Omega(t^n)$,
$\Omega^{n+1}:=\Omega(t^{n+1})$ and
$\partial\mathcal{E}^n := \cbr{(t,\boldsymbol{x})\in\partial\mathcal{E} :
  t\in I^n}$.

\begin{figure}[tbp]
\begin{center}
\includegraphics[width=5cm]{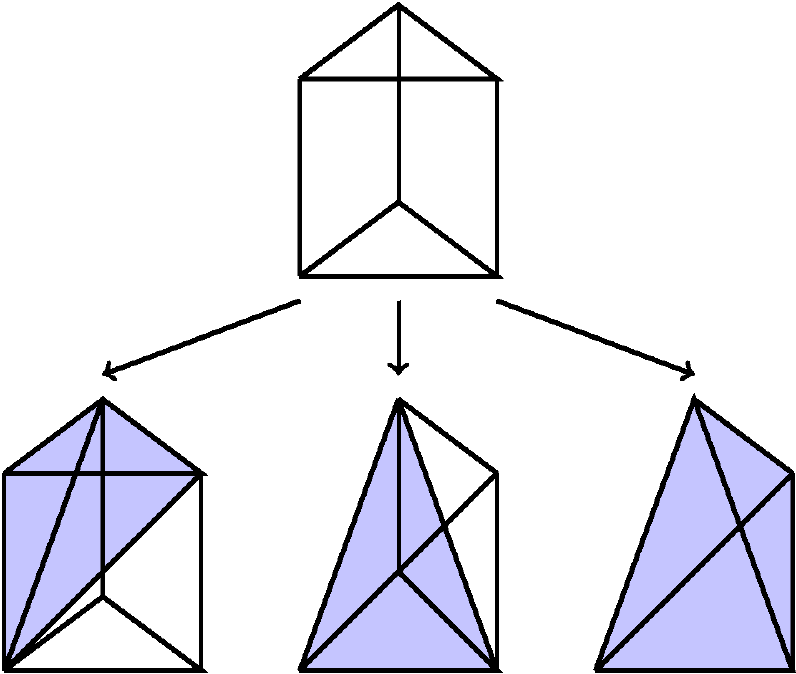}
\end{center}
\caption{The splitting of a space--time prism into three tetrahedra.\label{fig:prism2tet}}
\end{figure}

\begin{figure}[tbp]
\begin{center}
\includegraphics[width=10cm]{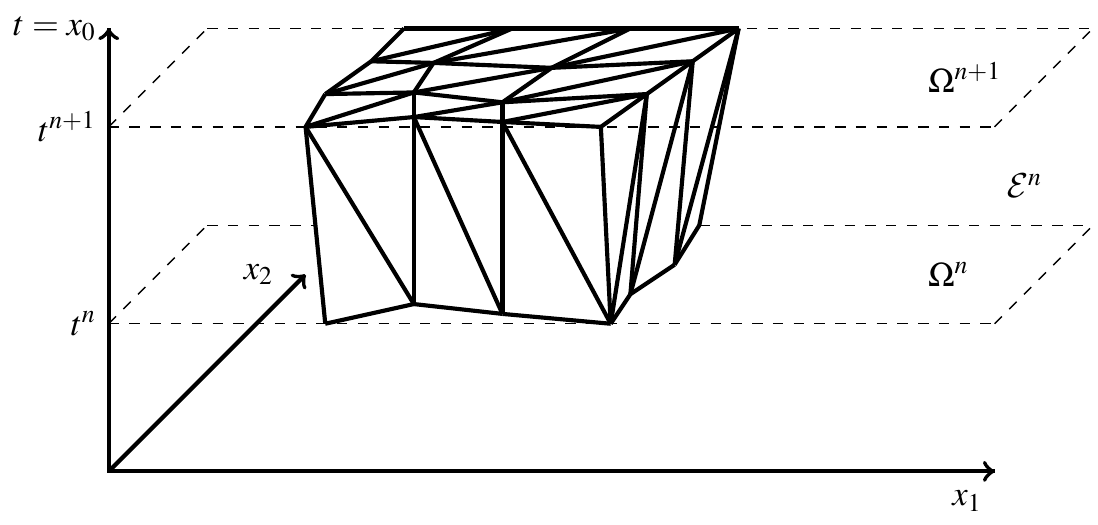}
\end{center}
\caption{A space--time tetrahedral mesh of a space--time slab
  $\mathcal{E}^n$.\label{fig:st-tet}}
\end{figure}

Let the evolution of the spatial domain during the time interval $I^n$
be represented by a sufficiently smooth and invertible mapping
$\Phi^n:\Omega(t^n) \to \Omega(t): \boldsymbol{x}\mapsto \Phi^n(\boldsymbol{x})$.
The standard approach~\cite{Klaij:2006, Masud:1997, Rhebergen:2012,
  Vegt:2002} to creating a space--time mesh in a space--time slab is
to extrude the spatial mesh of $\Omega^n$ to the new time level
$t=t^{n+1}$ according to the mapping $\Phi^n$. In the case of a
spatial simplicial mesh, this approach results in a mesh of the
space--time slab $\mathcal{E}^n$ consisting of space--time prisms. In
this paper, however, we follow the approach of~\cite{Ndri:2001,
  Ndri:2002, Wang:2015} and divide each space--time prism into three
space--time tetrahedra, see~\cref{fig:prism2tet}. The main advantage
of using space--time tetrahedra in this paper is the simplicity of
obtaining an approximate velocity field that is
$H({\rm div})$-conforming and point-wise divergence-free on
time-dependent domains. The triangulation of the space--time slab
$\mathcal{E}^n$ consisting of non-overlapping tetrahedral space--time
cells is denoted by $\mathcal{T}^n := \cbr{\mathcal{K}}$,
see~\cref{fig:st-tet}. The triangulation of the space--time domain
$\mathcal{E}$ is denoted by $\mathcal{T} := \cup_n\mathcal{T}^n$.

Consider now a single space--time cell
$\mathcal{K}_j \in \mathcal{T}^n$. The boundary of this space--time
cell is denoted by $\partial \mathcal{K}_j$. The outward unit
space--time normal vector on $\partial \mathcal{K}_j$ is given by
$(n_t^{\mathcal{K}_j}, \boldsymbol{n}^{\mathcal{K}_j})$. The boundary
$\partial\mathcal{K}_j$ may consist of a facet on which
$|n_t^{\mathcal{K}_j}| = 1$ (which we denote by $K^n_j$ if
$n_t^{\mathcal{K}_j} = -1$ or $K^{n+1}_j$ if
$n_t^{\mathcal{K}_j} = 1$) and
$\mathcal{Q}_{\mathcal{K}_j} = \partial \mathcal{K}_j \backslash
(K_j^n\cup K_j^{n+1})$.
In the remainder of this paper, we will drop the sub- and superscript
notatation when referring to the space--time normal vector and the
space--time cell wherever no confusion will occur.

In a space--time slab $\mathcal{E}^n$, the set of all facets for which
$\envert{n_t}\ne 1$ is denoted by $\mathcal{F}^n$ while the union of
these facets is denoted by $\Gamma^n$. The set $\mathcal{F}^n$ is
partitioned into a set of interior facets $\mathcal{F}_I^n$ and a set
of facets that lie on the boundary of $\mathcal{E}^n$,
$\mathcal{F}_B^n$, so that
$\mathcal{F}^n = \mathcal{F}_I^n \cup \mathcal{F}_B^n$. We furthermore
denote the set of facets that lie on the Neumann boundary,
$\partial\mathcal{E}^N \cap \partial\mathcal{E}^n$, by
$\mathcal{F}^n_N$.

We consider the following finite-dimensional function spaces on the
space--time slab $\mathcal{E}^n$:
\begin{subequations}
  \label{eq:st_general_func_spc}
  \begin{align}
  \label{eq:st_general_func_spc_a}
    \boldsymbol{V}_h^n &:= \cbr{ \boldsymbol{v}_h \in \sbr{L^2(\mathcal{T}^n)}^d, \
      \boldsymbol{v}_h \in \sbr{P_k(\mathcal{K})}^d\ \forall \mathcal{K} \in \mathcal{T}^n},
    \\
    \label{eq:st_general_func_spc_b}
    Q_h^n &:= \cbr{ q_h \in L^2(\mathcal{T}^n),\ q_h \in P_{k-1}(\mathcal{K})\
      \forall \mathcal{K} \in \mathcal{T}^n},
  \end{align}
\end{subequations}
where $P_l(D)$ denotes the space of polynomials of degree $l > 0$ on a
domain~$D$. The left and right traces of a function $q_h \in Q_h^n$ at
an interior facet $\mathcal{S} \in \mathcal{F}_I^n$ are denoted by
$q_h^l$ and $q_h^r$. In general $q_h^l \ne q_h^r$ so that it will be
useful to introduce the jump operator
$\jump{q_h\boldsymbol{n}} = q_h^l\boldsymbol{n}^l + q_h^r\boldsymbol{n}^r$. On a boundary
facet $\mathcal{S} \in \mathcal{F}_B^n$ the jump operator is defined
as $\jump{q_h\boldsymbol{n}} = q_h\boldsymbol{n}$. Similar expressions hold for
$\boldsymbol{v}_h \in \boldsymbol{V}_h^n$.

For the hybridized discontinuous Galerkin method, we require also
finite dimensional function spaces on $\Gamma^n$:
\begin{subequations}
  \label{eq:st_general_func_spc_facets}
  \begin{align}
  \label{eq:st_general_func_spc_facets_a}
    \bar{\boldsymbol{V}}_h^n &:= \cbr{ \bar{\boldsymbol{v}}_h \in \sbr{L^2(\mathcal{F}^n)}^d, \
      \bar{\boldsymbol{v}}_h \in \sbr{P_k(\mathcal{S})}^d
      \ \forall \mathcal{S} \in \mathcal{F}^n,\ \bar{\boldsymbol{v}}_h = 0 \ \text{on}\ \partial\mathcal{E}^D\cap\partial\mathcal{E}^n},
    \\
    \label{eq:st_general_func_spc_facets_d}
    \bar{Q}_h^n &:= \cbr{ \bar{q}_h \in L^2(\mathcal{F}^n),\
      \bar{q}_h \in P_{k}(\mathcal{S}) \ \forall \mathcal{S}\in\mathcal{F}^n}.
  \end{align}
\end{subequations}
For notational purposes, we introduce the spaces
$\boldsymbol{V}_h^{n,\star} = \boldsymbol{V}_h^n \times
\bar{\boldsymbol{V}}_h^n$ and
$Q_h^{n,\star} = Q_h^n \times \bar{Q}_h^n$. Functions pairs in
$\boldsymbol{V}_h^{n,\star}$ and $Q_h^{n,\star}$ will be denoted by
$\boldsymbol{v}_h^{\star} = (\boldsymbol{v}_h, \bar{\boldsymbol{v}}_h)
\in \boldsymbol{V}_h^{n,\star}$ and
$q_h^{\star} = (q_h, \bar{q}_h) \in Q_h^{n,\star}$.

\subsection{The finite element variational formulation}
\label{ss:variationalform}
We present now the finite element variational formulation of the
Navier--Stokes problem on time-dependent domains. For this, consider
first the steady Stokes problem on a time-dependent domain:
\begin{subequations}
  \label{eq:stokes}
  \begin{align}
    \label{eq:stokes_mom}
    - \nu \nabla^2 \boldsymbol{u} + \nabla p &= \boldsymbol{f}
    && \text{in} \ \mathcal{E},
    \\
    \label{eq:stokes_mass}
    \nabla \cdot \boldsymbol{u} &= 0 && \text{in} \ \mathcal{E},
  \end{align}
\end{subequations}
with boundary conditions
\begin{subequations}
  \label{eq:boundaryconditions_stokes}
  \begin{align}
    \label{eq:st_bc_dirichlet_stokes}
    \boldsymbol{u} &= \boldsymbol{0} && \text{on} \ \partial\mathcal{E}^D,
    \\
    \label{eq:st_bc_neumann_stokes}
    \del{p\mathbb{I} - \nu\nabla \boldsymbol{u}}\cdot\boldsymbol{n} &= \boldsymbol{h} && \text{on} \ \partial\mathcal{E}^N,
  \end{align}
\end{subequations}
with $\boldsymbol{h}:\partial\mathcal{E}^N\to\mathbb{R}^d$ given Neumann
boundary data. A straightforward extension of the hybridized
discontinuous Galerkin method of~\cite{Rhebergen:2017} to the Stokes
problem on time-dependent domains is given by: In each space--time
slab $\mathcal{E}^n$, $n=0,1,\cdots, N-1$, we seek
$(\boldsymbol{u}_h^{\star}, p_h^{\star}) \in \boldsymbol{V}_h^{n,\star} \times
Q_h^{n,\star}$ such that
\begin{equation}
  \label{eq:weakformulation_stokes}
    a_h^n(\boldsymbol{u}_h^{\star}, \boldsymbol{v}_h^{\star}) 
    + b_h^n(p_h^{\star}, \boldsymbol{v}_h^{\star}) - b_h^n(q_h^{\star}, \boldsymbol{u}_h^{\star}) =
  \sum_{\mathcal{K}\in\mathcal{T}^n} \int_{\mathcal{K}} \boldsymbol{f} \cdot \boldsymbol{v}_h \dif \boldsymbol{x}\dif t 
  - \sum_{\mathcal{S}\in\mathcal{F}^n_N}\int_{\mathcal{S}} \boldsymbol{h} \cdot \bar{\boldsymbol{v}}_h \dif s,
\end{equation}
for all
$(\boldsymbol{v}_h^{\star}, q_h^{\star}) \in \boldsymbol{V}_h^{n,\star} \times
Q_h^{n,\star}$, where
\begin{subequations}
  \begin{align}
    \label{eq:bilinform_a}
    a_h^n(\boldsymbol{u}^{\star}, \boldsymbol{v}^{\star})
    :=& 
        \sum_{\mathcal{K}\in\mathcal{T}^n} \int_{\mathcal{K}} \nu \nabla \boldsymbol{u} : \nabla \boldsymbol{v} \dif \boldsymbol{x} \dif t
        + \sum_{K\in\mathcal{T}^n} \int_{\mathcal{Q}_{\mathcal{K}}} 
        \frac{\nu \alpha}{h_{\mathcal{K}}}(\boldsymbol{u} - \bar{\boldsymbol{u}}) \cdot (\boldsymbol{v} - \bar{\boldsymbol{v}}) \dif s
    \\
    \nonumber
      &- \sum_{\mathcal{K}\in\mathcal{T}^n} \int_{\mathcal{Q}_{\mathcal{K}}} \nu \sbr{ (\boldsymbol{u} - \bar{\boldsymbol{u}}) 
        \cdot \pd{\boldsymbol{v}}{\boldsymbol{n}} 
        + \pd{\boldsymbol{u}}{\boldsymbol{n}} \cdot (\boldsymbol{v} - \bar{\boldsymbol{v}}) } \dif s,
    \\
    \label{eq:bilinform_b}
    b_h^n(p^{\star}, \boldsymbol{v}^{\star})
    := & 
         - \sum_{\mathcal{K}\in\mathcal{T}^n} \int_{\mathcal{K}} p \nabla\cdot \boldsymbol{v} \dif \boldsymbol{x} \dif t
         + \sum_{\mathcal{K}\in\mathcal{T}^n} \int_{\mathcal{Q}_{\mathcal{K}}} (\boldsymbol{v} - \bar{\boldsymbol{v}}) \cdot \boldsymbol{n}
         \bar{p} \dif s.
  \end{align}
\end{subequations}
Here $\alpha > 0$ is a penalty parameter that needs to be sufficiently
large to ensure stability. These bi-linear forms are similar to those
of~\cite{Rhebergen:2017, Labeur:2012, Rhebergen:2018a} with the
difference being that integration is over $d+1$-dimensional
space--time cells $\mathcal{K}$, as opposed to $d$-dimensional spatial
cells. Note furthermore that space--time slabs are completely
independent of each other due to the Stokes equations~\cref{eq:stokes}
not being time-dependent.

We will now discuss in more detail the variational formulation for the
convective parts of the linearized Navier--Stokes equations. In
particular, let $\boldsymbol{w}:\mathcal{E}\to\mathbb{R}^d$ be a given
divergence-free and $H(\text{div})$-conforming velocity field. We will
derive the discrete space--time variational formulation for
\begin{equation}
  \label{eq:navsto_convection}
  \partial_t\boldsymbol{u} + \nabla\cdot\del{\boldsymbol{u}\otimes\boldsymbol{w}} = \boldsymbol{0} 
  \qquad \text{in} \ \mathcal{E},
\end{equation}
with boundary condition
\begin{equation}
  \label{eq:st_bc_neumann_c}
  (n_t + \boldsymbol{w}\cdot\boldsymbol{n} )\boldsymbol{u} = \boldsymbol{r} \quad
  \text{on} \ \partial\mathcal{E}^{-},  
\end{equation}
where $\partial\mathcal{E}^{-}$ is the part of $\partial\mathcal{E}^N$
where $n_t + \boldsymbol{w}\cdot\boldsymbol{n} < 0$, and where
$\boldsymbol{r}:\partial\mathcal{E}^{-}\to\mathbb{R}^d$ is given boundary
data. In each space--time slab $\mathcal{E}^n$, we
multiply~\cref{eq:navsto_convection} by a test function
$\boldsymbol{v}_h \in \boldsymbol{V}_h^n$, integrate over a cell
$\mathcal{K}\in\mathcal{T}^n$, approximate $\boldsymbol{w}$ by a
divergence-free $\boldsymbol{w}_h \in \boldsymbol{V}_h^n \cap H({\rm div})$, apply
Green's identity in space--time and sum over all cells of the
triangulation and over all space--time slabs,
\begin{equation}
  \label{eq:firststep_conv}
  \sum_{n=0}^{N-1}\bigg(
  -\sum_{\mathcal{K}\in\mathcal{T}^n}\int_{\mathcal{K}}\del{\boldsymbol{u}_h\cdot\partial_t\boldsymbol{v}_h + \boldsymbol{u}_h\otimes \boldsymbol{w}_h : \nabla \boldsymbol{v}_h}
  \dif \boldsymbol{x}\dif t 
+ \sum_{\mathcal{K}\in\mathcal{T}^n}\int_{\partial\mathcal{K}} H_{\mathcal{K}}^n(\boldsymbol{u}_h^{\star}, \boldsymbol{w}_h; n_t, \boldsymbol{n}) \cdot \boldsymbol{v}_h \dif s\bigg)
  = 0.
\end{equation}
For stability purposes, the convective flux on the cell boundary
$\partial\mathcal{K}$ in the space--time normal direction,
$\del{n_t + \boldsymbol{w}_h\cdot\boldsymbol{n}}\boldsymbol{u}_h$, was replaced by the
space--time upwind flux
\begin{equation}
  \label{eq:numericalflux}
  H_{\mathcal{K}}^n(\boldsymbol{u}_h^{\star}, \boldsymbol{w}_h; n_t, \boldsymbol{n})
  =
  \begin{cases}
    \del{n_t + \boldsymbol{w}_h\cdot\boldsymbol{n}}\del{\boldsymbol{u}_h + \lambda(\bar{\boldsymbol{u}}_h - \boldsymbol{u}_h)} & \text{on} \ \mathcal{Q}_{\mathcal{K}},
    \\
    \boldsymbol{u}_hn_t   & \text{on}\ K^{n+1},
    \\
    \boldsymbol{u}_h^-n_t & \text{on}\ K^n,
  \end{cases}
\end{equation}
where, on $K^n$,
$\boldsymbol{u}_h^- = \lim_{\varepsilon\to 0} \boldsymbol{u}_h(t_n - \varepsilon)$,
and where $\lambda = 1$ if $n_t + \boldsymbol{w}_h\cdot\boldsymbol{n} < 0$ and
$\lambda = 0$ otherwise. Substituting this expression
into~\cref{eq:firststep_conv},
\begin{multline}
  \label{eq:firststep_conv_flux}
  \sum_{n=0}^{N-1}\bigg(
  -\sum_{\mathcal{K}\in\mathcal{T}^n}\int_{\mathcal{K}}\del{\boldsymbol{u}_h\cdot\partial_t\boldsymbol{v}_h + \boldsymbol{u}_h\otimes \boldsymbol{w}_h : \nabla \boldsymbol{v}_h}
  \dif \boldsymbol{x}\dif t 
  + \sum_{\mathcal{K}\in\mathcal{T}^n}\int_{K^{n+1}}\boldsymbol{u}_h\cdot\boldsymbol{v}_h \dif \boldsymbol{x} 
  - \sum_{\mathcal{K}\in\mathcal{T}^n}\int_{K^{n}}\boldsymbol{u}_h^-\cdot\boldsymbol{v}_h \dif \boldsymbol{x}
  \\ 
  + \sum_{\mathcal{K}\in\mathcal{T}^n}\int_{\mathcal{Q}_{\mathcal{K}}} 
  H_{\mathcal{K}}^n(\boldsymbol{u}_h^{\star}, \boldsymbol{w}_h; n_t, \boldsymbol{n}) \cdot \boldsymbol{v}_h \dif s\bigg)
  = 0,
\end{multline}
where we used that $n_t=1$ on $K^{n+1}$ and $n_t=-1$ on $K^n$. Note
that the numerical flux on the boundary of cell $\mathcal{K}$ depends
only on the local cell unknown $\boldsymbol{u}_h$ and the facet
unknown $\bar{\boldsymbol{u}}_h$. As such, the numerical flux on an
interior facet $\mathcal{S} \in \mathcal{F}$, shared by two adjacent
cells $\mathcal{K}^+$ and $\mathcal{K}^-$, need not be equal, i.e.,
$H_{\mathcal{K}^+}^n((\boldsymbol{u}_h^+, \bar{\boldsymbol{u}}_h),
\boldsymbol{w}_h^+; n_t^+, \boldsymbol{n}^+) \ne
H_{\mathcal{K}^-}^n((\boldsymbol{u}_h^-, \bar{\boldsymbol{u}}_h),
\boldsymbol{w}_h^-; n_t^-, \boldsymbol{n}^-)$.  In order to guarantee
local conservation, we follow~\cite{Nguyen:2011} and impose that the
$L^2$-projection of the space--time normal component of the numerical
flux into $\bar{\boldsymbol{V}}_h^n$ is single-valued:
\begin{equation}
  \label{eq:singlevaluedflux}
  \sum_{n=0}^{N-1}\del{
  \int_{\mathcal{K}\in\mathcal{T}^n}H_{\mathcal{K}}^n(\boldsymbol{u}_h^{\star}, \boldsymbol{w}_h; n_t, \boldsymbol{n}) \cdot \bar{\boldsymbol{v}}_h \dif s
  - \int_{\partial\mathcal{E}^+} (n_t + \bar{\boldsymbol{w}}_h\cdot \boldsymbol{n}) \bar{\boldsymbol{u}}_h\cdot \bar{\boldsymbol{v}}_h \dif s} =
 \sum_{n=0}^{N-1}\int_{\partial\mathcal{E}^-} \boldsymbol{r} \cdot \bar{\boldsymbol{v}}_h \dif s,
\end{equation}
where $\partial\mathcal{E}^+$ is the part of $\partial\mathcal{E}^N$
where $n_t + \boldsymbol{w}\cdot\boldsymbol{n} \ge 0$ and where we used the inflow
boundary condition~\cref{eq:st_bc_neumann_c}.
Subtracting~\cref{eq:singlevaluedflux}
from~\cref{eq:firststep_conv_flux}, and noting that each space--time
slab depends only on the previous space--time slab, so that we may
drop the summation over space--time slabs, we find the following
finite element variational formulation for~\cref{eq:navsto_convection}
and~\cref{eq:st_bc_neumann_c}: In each space--time slab
$\mathcal{E}^n$, $n=0, \cdots, N-1$, we seek
$\boldsymbol{u}_h^{\star} \in \boldsymbol{V}_h^{n,\star}$ such that
\begin{equation}
  \label{eq:weakformconvection}
  t_h^n(\boldsymbol{u}_h^{\star}, \boldsymbol{w}_h^{\star}, \boldsymbol{v}_h^{\star}) 
  = \sum_{\mathcal{K}\in\mathcal{T}^n}\int_{K^{n}}\boldsymbol{u}_h^-\cdot\boldsymbol{v}_h \dif \boldsymbol{x} 
  - \int_{\partial\mathcal{E}^-} \boldsymbol{r} \cdot \bar{\boldsymbol{v}}_h \dif s,
\end{equation}
where the tri-linear form for the convective term is defined as:
\begin{multline}
  \label{eq:trilinearform}
  t_h^n(\boldsymbol{u}^{\star}, \boldsymbol{w}^{\star}, \boldsymbol{v}^{\star}) = 
  -\sum_{\mathcal{K}\in\mathcal{T}^n}\int_{\mathcal{K}}\del{\boldsymbol{u}\cdot\partial_t\boldsymbol{v} + \boldsymbol{u}\otimes \boldsymbol{w} : \nabla \boldsymbol{v}}
  \dif \boldsymbol{x}\dif t 
  + \sum_{\mathcal{K}\in\mathcal{T}^n}\int_{K^{n+1}}\boldsymbol{u}\cdot\boldsymbol{v} \dif \boldsymbol{x}
  \\
  + \sum_{\mathcal{K}\in\mathcal{T}^n}\int_{\mathcal{Q}_{\mathcal{K}}} 
  H_{\mathcal{K}}^n(\boldsymbol{u}^{\star}, \boldsymbol{w}; n_t, \boldsymbol{n}) \cdot \del{\boldsymbol{v} - \bar{\boldsymbol{v}}} \dif s
  + \int_{\partial\mathcal{E}^+} (n_t + \bar{\boldsymbol{w}}\cdot \boldsymbol{n}) \bar{\boldsymbol{u}}\cdot \bar{\boldsymbol{v}} \dif s.
\end{multline}
We remark that for $n=0$, $\boldsymbol{u}_h^-$ is the projection of the
initial condition~\cref{eq:st_ic} into
$\boldsymbol{V}_h^0 \cap H(\text{div})$ such that $\boldsymbol{u}_h^-$ is point-wise
divergence free.

Combining~\cref{eq:weakformulation_stokes}
and~\cref{eq:weakformconvection}, we conclude this section by stating
the discontinuous Galerkin in time and space--time hybridized
discontinuous Galerkin variational formulation for the incompressible
Navier--Stokes problem~\cref{eq:navsto}--\cref{eq:boundaryconditions}:
In each space--time slab $\mathcal{E}^n$, $n=0,\cdots, N-1$, we seek
$(\boldsymbol{u}_h^{\star}, p_h^{\star}) \in \boldsymbol{V}_h^{n,\star} \times
Q_h^{n,\star}$ such that
\begin{equation}
  \label{eq:weakformulation_stns}
  t_h^n(\boldsymbol{u}_h^{\star}, \boldsymbol{u}_h^{\star}, \boldsymbol{v}_h^{\star}) + 
    a_h^n(\boldsymbol{u}_h^{\star}, \boldsymbol{v}_h^{\star}) + b_h^n(p_h^{\star}, \boldsymbol{v}_h^{\star}) - b_h^n(q_h^{\star}, \boldsymbol{u}_h^{\star}) =
  \sum_{\mathcal{K}\in\mathcal{T}^n} \int_{\mathcal{K}} \boldsymbol{f} \cdot \boldsymbol{v}_h \dif \boldsymbol{x} \dif t
  - \sum_{\mathcal{S}\in\mathcal{F}^n_N}\int_{\mathcal{S}} \boldsymbol{g} \cdot \bar{\boldsymbol{v}}_h \dif s 
  + \int_{\Omega_n}\boldsymbol{u}_h^- \cdot \boldsymbol{v}_h \dif s,
\end{equation}
for all
$(\boldsymbol{v}_h^{\star}, q_h^{\star}) \in \boldsymbol{V}_h^{n,\star} \times
Q_h^{n,\star}$.

\paragraph{Arbitrary Lagrangian Eulerian formulation} For practical
implementations it may be preferable to consider the space--time
hybridized discontinuous Galerkin variational
formulation~\cref{eq:weakformulation_stns} in arbitrary Lagrangian
Eulerian (ALE) formulation. The ALE formulation is obtained by noting
that the space--time normal on space--time cell boundaries
$\mathcal{Q}_{\mathcal{K}}$ may be written as
$(n_t, \boldsymbol{n}) = (-\boldsymbol{v}_g\cdot \boldsymbol{n}, \boldsymbol{n})$, where
$\boldsymbol{v}_g$ is the grid velocity~\cite{Vegt:2002}. Only the tri-linear
form, $t_h^n(\boldsymbol{u}^{\star}, \boldsymbol{w}^{\star}, \boldsymbol{v}^{\star})$,
needs to be rewritten in ALE formulation, since the other terms in the
variational formulation do not depend on $n_t$. The ALE formulation of
$t_h^n(\boldsymbol{u}^{\star}, \boldsymbol{w}^{\star}, \boldsymbol{v}^{\star})$ is given by
\begin{multline}
  \label{eq:trilinearform_ale}
  t_h^n(\boldsymbol{u}^{\star}, \boldsymbol{w}^{\star}, \boldsymbol{v}^{\star}) = 
  -\sum_{\mathcal{K}\in\mathcal{T}^n}\int_{\mathcal{K}}\del{\boldsymbol{u}\cdot\partial_t\boldsymbol{v} + \boldsymbol{u}\otimes \boldsymbol{w} : \nabla \boldsymbol{v}}
  \dif \boldsymbol{x}\dif t 
  + \sum_{\mathcal{K}\in\mathcal{T}^n}\int_{K^{n+1}}\boldsymbol{u}\cdot\boldsymbol{v} \dif \boldsymbol{x}
  \\
  + \sum_{\mathcal{K}\in\mathcal{T}^n}\int_{\mathcal{Q}_{\mathcal{K}}} 
  G_{\mathcal{K}}^n(\boldsymbol{u}^{\star}, \boldsymbol{w}; \boldsymbol{v}_g, \boldsymbol{n}) \cdot \del{\boldsymbol{v} - \bar{\boldsymbol{v}}} \dif s
  + \int_{\partial\mathcal{E}^+} \del{(\bar{\boldsymbol{w}}-\boldsymbol{v}_g)\cdot \boldsymbol{n}} \bar{\boldsymbol{u}}\cdot \bar{\boldsymbol{v}} \dif s,
\end{multline}
with numerical flux
\begin{equation}
  \label{eq:numericalflux_G}
  G_{\mathcal{K}}^n(\boldsymbol{u}^{\star}, \boldsymbol{w}; \boldsymbol{v}_g, \boldsymbol{n})
  =
  \del{\boldsymbol{u}_h + \lambda(\bar{\boldsymbol{u}}_h - \boldsymbol{u}_h)} (\boldsymbol{w}_h-\boldsymbol{v}_g)\cdot\boldsymbol{n}.
\end{equation}

\section{Properties of the discrete variational formulation}
\label{s:properties}

To find the solution
$(\boldsymbol{u}_h^{\star}, p_h^{\star}) \in \boldsymbol{V}_h^{n,\star} \times
Q_h^{n,\star}$
to the non-linear variational
formulation~\cref{eq:weakformulation_stns}, we use a Picard iteration
scheme: in every space--time slab, given a solution
$(\boldsymbol{u}_h^{\star,k}, p_h^{\star,k})$ we seek a solution
$(\boldsymbol{u}_h^{\star,k+1}, p_h^{\star,k+1})$ in iteration $k+1$ that
solves the linear variational formulation
\begin{multline}
  \label{eq:linearweakformulation_stns}
  t_h^n(\boldsymbol{u}_h^{\star,k+1}, \boldsymbol{u}_h^{\star,k}, \boldsymbol{v}_h^{\star}) + 
    a_h^n(\boldsymbol{u}_h^{\star,k+1}, \boldsymbol{v}_h^{\star}) + b_h^n(p_h^{\star,k+1}, \boldsymbol{v}_h^{\star}) - b_h^n(q_h^{\star}, \boldsymbol{u}_h^{\star,k+1}) =
  \\
  \sum_{\mathcal{K}\in\mathcal{T}^n} \int_{\mathcal{K}} \boldsymbol{f} \cdot \boldsymbol{v}_h \dif \boldsymbol{x} \dif t
  - \sum_{\mathcal{S}\in\mathcal{F}^n_N}\int_{\mathcal{S}} \boldsymbol{g} \cdot \bar{\boldsymbol{v}}_h \dif s 
  + \int_{\Omega_n}\boldsymbol{u}_h^- \cdot \boldsymbol{v}_h \dif s,
\end{multline}
for all
$(\boldsymbol{v}_h^{\star}, q_h^{\star}) \in \boldsymbol{V}_h^{n,\star} \times
Q_h^{n,\star}$.
The iterations are stopped when a certain convergence criterium has
been met, at which point we set
$(\boldsymbol{u}_h^{\star}, p_h^{\star}) = (\boldsymbol{u}_h^{\star,k+1},
p_h^{\star,k+1})$.

The approximate velocity field $\boldsymbol{u}_h^k$ that is obtained at each
Picard iteration $k$ is $H({\rm div})$-conforming and point-wise
divergence-free in each space--time cell
$\mathcal{K}\in\mathcal{T}^n$. To see this, we note that by taking
$(\boldsymbol{v}_h, \bar{\boldsymbol{v}}_h, \bar{q}_h) = (\boldsymbol{0}, \boldsymbol{0}, 0)$ and
$q_h = \nabla\cdot\boldsymbol{u}_h^k \in Q_h^n$
in~\cref{eq:weakformulation_stns} results in
\begin{equation}
  \int_{\mathcal{K}} (\nabla\cdot\boldsymbol{u}_h^k)^2 \dif \boldsymbol{x} \dif t = 0 \qquad 
  \forall\mathcal{K}\in\mathcal{T}^n,
\end{equation}
from which it follows that $\nabla\cdot\boldsymbol{u}_h^k=0$ for
$(\boldsymbol{x},t)\in\mathcal{K}$ and for each
$\mathcal{K}\in\mathcal{T}^n$. Furthermore, taking
$(\boldsymbol{v}_h, \bar{\boldsymbol{v}}_h, q_h) = (\boldsymbol{0}, \boldsymbol{0}, 0)$ and
$\bar{q}_h = \jump{(\boldsymbol{u}_h^k -
  \bar{\boldsymbol{u}}_h^k)\cdot\boldsymbol{n}} \in \bar{Q}_h^n$
in~\cref{eq:weakformulation_stns} results in
\begin{equation}
  \sum_{\mathcal{S}\in\mathcal{F}_I^n}\int_{\mathcal{S}} \jump{\boldsymbol{u}_h^k\cdot\boldsymbol{n}}^2 \dif \boldsymbol{x} \dif t
  +
  \sum_{\mathcal{S}\in\mathcal{F}_B^n}\int_{\mathcal{S}} \del{ ( \boldsymbol{u}_h^k - \bar{\boldsymbol{u}}_h^k )\cdot\boldsymbol{n} }^2 \dif \boldsymbol{x} \dif t
  = 0,
\end{equation}
where we used that $\bar{\boldsymbol{u}}_h^k$ is single-valued on interior
facets. It follows that $\boldsymbol{u}_h^k \cdot \boldsymbol{n}$ is single-valued
on interior facets and
$\boldsymbol{u}_h^k\cdot\boldsymbol{n} = \bar{\boldsymbol{u}}_h^k\cdot \boldsymbol{n}$ on boundary
facets, i.e., $\boldsymbol{u}_h^k$ is $H({\rm div})$-conforming.

Many finite element methods result in discretizations of the
Navier--Stokes equations that are either energy stable or locally
conservative. The space--time HDG
method~\cref{eq:weakformulation_stns}, however, is both
simultaneously, even on time-dependent domains. To see this, note that
by~\cref{eq:singlevaluedflux} the $L^2$-projection of the space--time
normal component of the numerical flux into $\bar{\boldsymbol{V}}_h^n$ is
single-valued guaranteeing that the variational
formulation~\cref{eq:weakformulation_stns} is locally conservative. We
next show that the space--time variational
formulation~\cref{eq:weakformulation_stns} is also energy stable.

Consider~\cref{eq:linearweakformulation_stns} in the first space--time
slab $\mathcal{E}_n$, $n=0$. Assume we are given a solution
$(\boldsymbol{u}_h^k, \bar{\boldsymbol{u}}_h^k, p_h^k, \bar{p}_h^k)$ from Picard
iteration $k$. For notational reasons, set
$\boldsymbol{u}_h^{\star} = \boldsymbol{u}_h^{\star,k+1}$,
$\boldsymbol{w}_h^{\star} = \boldsymbol{u}_h^{\star,k}$ and
$p_h^{\star} = p_h^{\star,k+1}$. For homogeneous boundary conditions,
$f=0$, and taking
$(\boldsymbol{v}_h^{\star}, q_h^{\star}) = (\boldsymbol{u}_h^{\star}, p_h^{\star})$
in~\cref{eq:linearweakformulation_stns},
\begin{equation}
  \label{eq:weakformulation_stns_1}
  t_h^0(\boldsymbol{u}_h^{\star}, \boldsymbol{w}_h^{\star}, \boldsymbol{u}_h^{\star}) + 
    a_h^0(\boldsymbol{u}_h^{\star}, \boldsymbol{u}_h^{\star})  =
    \int_{\Omega_0} \boldsymbol{u}_h^- \cdot \boldsymbol{u}_h \dif s,
\end{equation}
with the projection of the initial condition~\cref{eq:st_ic} into
$\boldsymbol{V}_h^0 \cap H(\text{div})$ is such that $\boldsymbol{u}_h^-$ is
point-wise divergence-free. Consider first the tri-linear form,
$t_h^0(\boldsymbol{u}_h^{\star}, \boldsymbol{w}_h^{\star}, \boldsymbol{u}_h^{\star})$, and
note that
\begin{multline}
  \label{eq:trilin_uuu}
  t_h^0(\boldsymbol{u}_h^{\star}, \boldsymbol{w}_h^{\star}, \boldsymbol{u}_h^{\star}) := 
  \sum_{\mathcal{K}\in\mathcal{T}^0} \int_{K^{1}}\envert{\boldsymbol{u}_h}^2 \dif \boldsymbol{x}
  - \sum_{\mathcal{K}\in\mathcal{T}^0} \int_{\mathcal{K}} \tfrac{1}{2}\del{\partial_t\envert{\boldsymbol{u}_h}^2 
    + \nabla \cdot \del{(\boldsymbol{u}_h\otimes \boldsymbol{w}_h)\cdot \boldsymbol{u}_h}} \dif \boldsymbol{x} \dif t
  \\
  + \sum_{\mathcal{K}\in\mathcal{T}^0} \int_{\mathcal{Q}^0_{\mathcal{K}}} \tfrac{1}{2} (n_t + \boldsymbol{w}_h\cdot \boldsymbol{n}) \envert{\boldsymbol{u}_h}^2 \dif s
  - \sum_{\mathcal{K}\in\mathcal{T}^0} \int_{\mathcal{Q}^0_{\mathcal{K}}} \tfrac{1}{2} (n_t + \boldsymbol{w}_h\cdot \boldsymbol{n}) \envert{\bar{\boldsymbol{u}}_h}^2 \dif s
  \\
  + \int_{\partial\mathcal{E}^N \cap I_0}\max(n_t + \bar{\boldsymbol{w}}_h\cdot \boldsymbol{n}, 0)\envert{\bar{\boldsymbol{u}}_h}^2 \dif s
  + \sum_{\mathcal{K}\in\mathcal{T}^0} \int_{\mathcal{Q}^0_{\mathcal{K}}} \tfrac{1}{2}
  \envert{n_t + \boldsymbol{w}_h\cdot \boldsymbol{n}} \envert{\boldsymbol{u}_h - \bar{\boldsymbol{u}}_h}^2 \dif s.
\end{multline}
Using that
\begin{equation}
  \label{eq:intbypartsterm_t}
  \sum_{\mathcal{K}\in\mathcal{T}^0} \int_{\mathcal{K}} \partial_t\envert{\boldsymbol{u}_h}^2 \dif \boldsymbol{x}\dif t
  =
  \sum_{\mathcal{K}\in\mathcal{T}^0} \int_{K^1} \envert{\boldsymbol{u}_h}^2 \dif \boldsymbol{x} 
  - \sum_{\mathcal{K}\in\mathcal{T}^0} \int_{K^0} \envert{\boldsymbol{u}_h}^2 \dif \boldsymbol{x}
  + \sum_{\mathcal{K}\in\mathcal{T}^0} \int_{\mathcal{Q}^0_{\mathcal{K}}} \envert{\boldsymbol{u}_h}^2 n_t \dif s,  
\end{equation}
and
\begin{align}
  \label{eq:intbypartsterm_x}
  \sum_{\mathcal{K}\in\mathcal{T}^0} \int_{\mathcal{K}} 
  \nabla \cdot \del{(\boldsymbol{u}_h\otimes \boldsymbol{w}_h)\cdot \boldsymbol{u}_h} \dif x
  &=
    \sum_{\mathcal{K}\in\mathcal{T}^0}\int_{\mathcal{Q}^0_{\mathcal{K}}}\boldsymbol{w}_h\cdot \boldsymbol{n}\envert{\boldsymbol{u}_h}^2\dif s,
  \\    
  \label{eq:intbypartsterm_3}
  \sum_{\mathcal{K}\in\mathcal{T}^0} \int_{\mathcal{Q}^0_{\mathcal{K}}} 
  (n_t + \boldsymbol{w}_h\cdot \boldsymbol{n}) \envert{\bar{\boldsymbol{u}}_h}^2 \dif s    
  &= \int_{\partial\mathcal{E}^N\cap I_0}
    (n_t + \bar{\boldsymbol{w}}_h\cdot \boldsymbol{n}) \envert{\bar{\boldsymbol{u}}_h}^2 \dif s,
\end{align}
where the second equality is by single-valuedness of $\bar{\boldsymbol{u}}_h$
on facets and since $\boldsymbol{w}_h$ is $H({\rm
  div})$-conforming.
Combining~\cref{eq:trilin_uuu}--\cref{eq:intbypartsterm_3}, and
simplifying terms,
\begin{multline}
  \label{eq:trilin_uuu_comb}
  t_h^0(\boldsymbol{u}_h^{\star}, \boldsymbol{w}_h^{\star}, \boldsymbol{u}_h^{\star}) =
  \sum_{\mathcal{K}\in\mathcal{T}^0} \int_{K^1} \tfrac{1}{2} \envert{\boldsymbol{u}_h}^2 \dif \boldsymbol{x}
  + 
  \sum_{\mathcal{K}\in\mathcal{T}^0} \int_{K^0} \tfrac{1}{2} \envert{\boldsymbol{u}_h}^2 \dif \boldsymbol{x}
  \\
  +
  \sum_{\mathcal{K}\in\mathcal{T}^0} \int_{\mathcal{Q}^0_{\mathcal{K}}} \tfrac{1}{2} 
  \envert{n_t + \boldsymbol{w}_h\cdot \boldsymbol{n}} \envert{\boldsymbol{u}_h - \bar{\boldsymbol{u}}_h}^2 \dif s
  + \int_{\partial\mathcal{E}^N\cap I_0}\frac{1}{2}
  \envert{n_t + \bar{\boldsymbol{w}}_h\cdot \boldsymbol{n}} \envert{\bar{\boldsymbol{u}}_h}^2 \dif s.
\end{multline}
Combining~\cref{eq:trilin_uuu_comb}
with~\cref{eq:weakformulation_stns_1}, and using that
$a_h^0(\boldsymbol{u}_h^{\star}, \boldsymbol{u}_h^{\star}) \ge 0$ for $\alpha$ large
enough (see~\cite{Labeur:2012, Rhebergen:2018a, Rhebergen:2017}), we
obtain
\begin{equation}
  \label{eq:almostthere}
  \int_{\Omega_1} \tfrac{1}{2} \envert{\boldsymbol{u}_h}^2 \dif \boldsymbol{x}
  + 
  \int_{\Omega_0} \tfrac{1}{2} \envert{\boldsymbol{u}_h}^2 \dif \boldsymbol{x}
  -
  \int_{\Omega_0} \boldsymbol{u}_h^- \cdot \boldsymbol{u}_h \dif \boldsymbol{x} \le 0.
\end{equation}
Using that
$2\boldsymbol{u}_h^- \cdot \boldsymbol{u}_h = |\boldsymbol{u}_h|^2 +
|\boldsymbol{u}_h^-|^2 -|\boldsymbol{u}_h - \boldsymbol{u}_h^-|^2$,
energy stability follows for each Picard iteration:
\begin{equation}
  \label{eq:energystab2}
  \int_{\Omega_1} \envert{\boldsymbol{u}_h}^2 \dif \boldsymbol{x}
  \le \int_{\Omega_0} \envert{\boldsymbol{u}_h^-}^2 \dif \boldsymbol{x}.
\end{equation}
Energy stability is now proven by induction for all $n > 0$ by using
$\boldsymbol{u}_h$ from space--time slab $\mathcal{E}_{n-1}$ as initial
condition for the variational
formulation~\cref{eq:linearweakformulation_stns} in space--time slab
$\mathcal{E}_n$. We have therefore shown that our space--time
variational formulation is both energy stable \emph{and} locally
conservative, even on dynamic meshes.

\section{Numerical examples}
\label{s:numericalexamples}
All simulations were implemented using the Modular Finite Element
Method (MFEM) library~\cite{mfem-library}. Furthermore, for all
simulations we use the penalty parameter $\alpha = 6k^2$.

In each space--time slab $\mathcal{E}^n$ ($n=0,\cdots, N-1$) we solve
the Navier--Stokes equations by Picard
iteration~\cref{eq:linearweakformulation_stns} with stopping criterion
\begin{equation}
  \label{eq:stopping_criterion}
  \max \left\{\frac{\|\boldsymbol{u}^k_h - \boldsymbol{u}^{k-1}_h\|_\infty}{\|\boldsymbol{u}^k_h - \boldsymbol{u}^0_h\|_\infty}, 
    \frac{\|p^k_h - p^{k-1}_h\|_\infty}{\|p^k_h - p^0_h\|_\infty} \right\} < \text{TOL},
\end{equation}
where $\|\cdot \|_{\infty}$ is the discrete $l^\infty$-norm and
$\text{TOL}$ is a user given parameter. 

Let
$U\in \mathbb{R}^{\dim \boldsymbol{V}^n_h}, P\in \mathbb{R}^{\dim
  Q^n_h}, \bar{U}\in \mathbb{R}^{\dim \bar{\boldsymbol{V}}^n_h},
\bar{P}\in \mathbb{R}^{\dim \bar{Q}^n_h}$ be the vectors of
coefficients of
${\boldsymbol{u}}_h, p_h, \bar{\boldsymbol{u}}_h, \bar{p}_h$ with
respect to the basis of the corresponding vector spaces. Then
$W^T = [U^T\ P^T]$ is the vector of all element degrees-of-freedom and
$\bar{W}^T = [\bar{U}^T\ \bar{P}^T]$ is the vector of all facet
degrees-of-freedom. At each Picard
iteration~\cref{eq:linearweakformulation_stns} a linear system of the
following block-matrix structure needs to be solved:
\begin{equation}
  \label{eq:before_SC}
  \begin{bmatrix}
    A & B \\ C & D
  \end{bmatrix}
  \begin{bmatrix}
    W \\ \bar{W}
  \end{bmatrix}
  =
  \begin{bmatrix}
    F \\ \bar{F}
  \end{bmatrix}.
\end{equation}
As with all other hybridizable discontinuous Galerkin methods, $A$ has
a block-diagonal structure. It is therefore cheap to eliminate $W$
from~\cref{eq:before_SC} to obtain the reduced linear system
$(-CA^{-1}B + D) \bar{W} = \bar{F} - CA^{-1}F$. We use the direct
solver of MUMPS~\cite{MUMPS:1,MUMPS:2} through
PETSc~\cite{petsc-web-page,petsc-user-ref,petsc-efficient} to solve
this system of linear equations. Given $\bar{W}$ we can then compute
$W$ cell-wise according to $W = A^{-1}(F - B\bar{W})$.

\subsection{Convergence rates and pressure-robustness}
\label{ss:convrates}
\begin{figure}[tbp]
  \begin{center}
    \includegraphics[width=.3\linewidth]{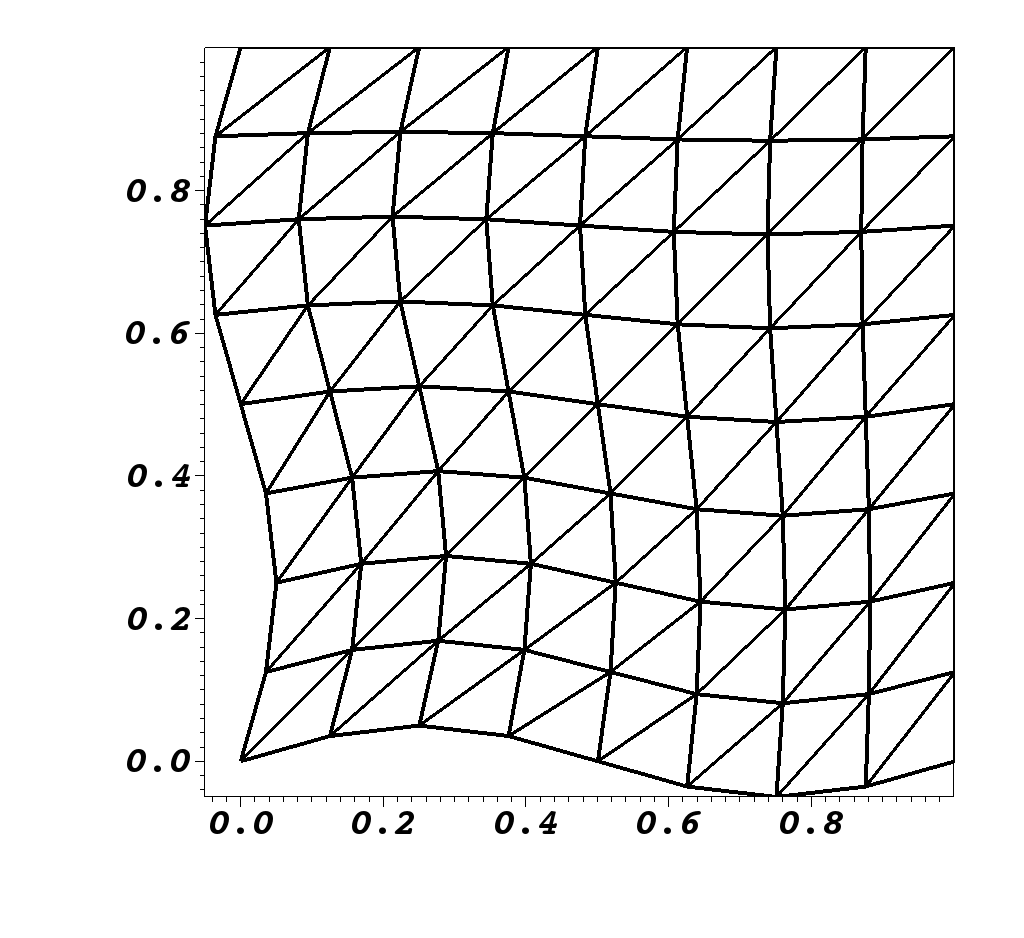}
    \includegraphics[width=.3\linewidth]{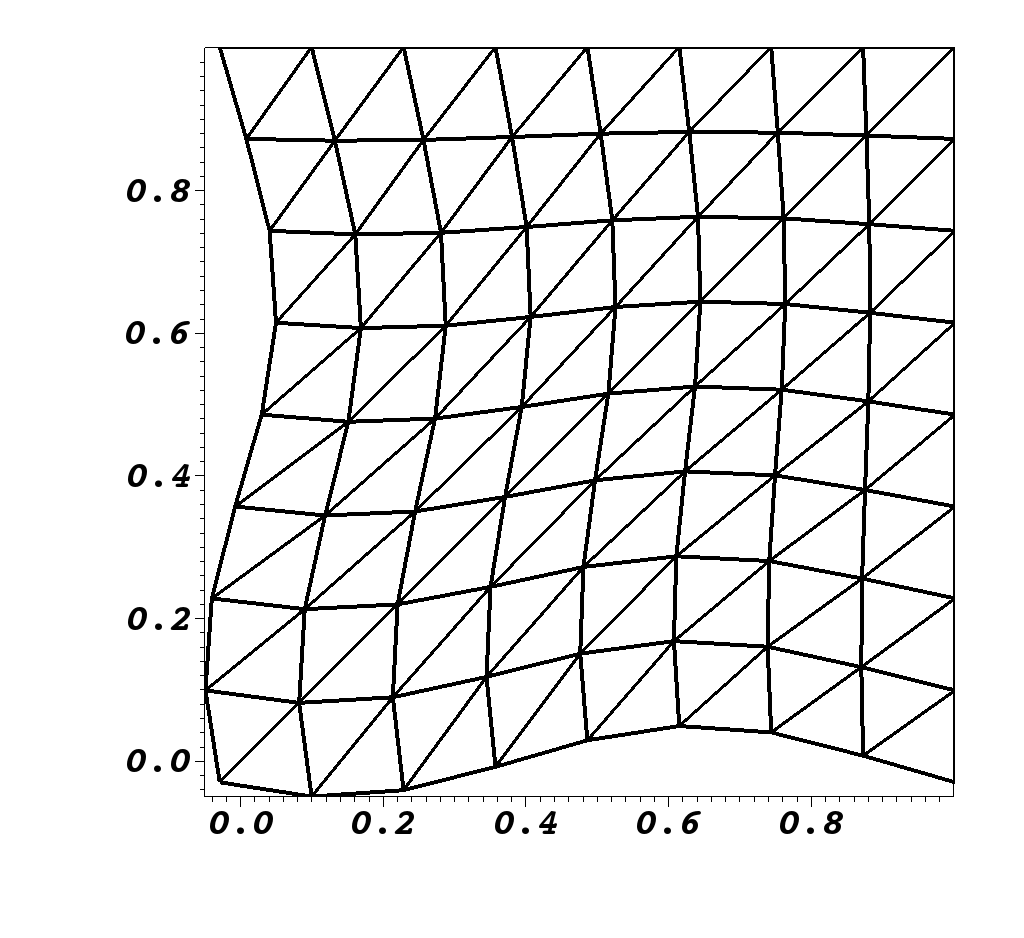}
    \includegraphics[width=.3\linewidth]{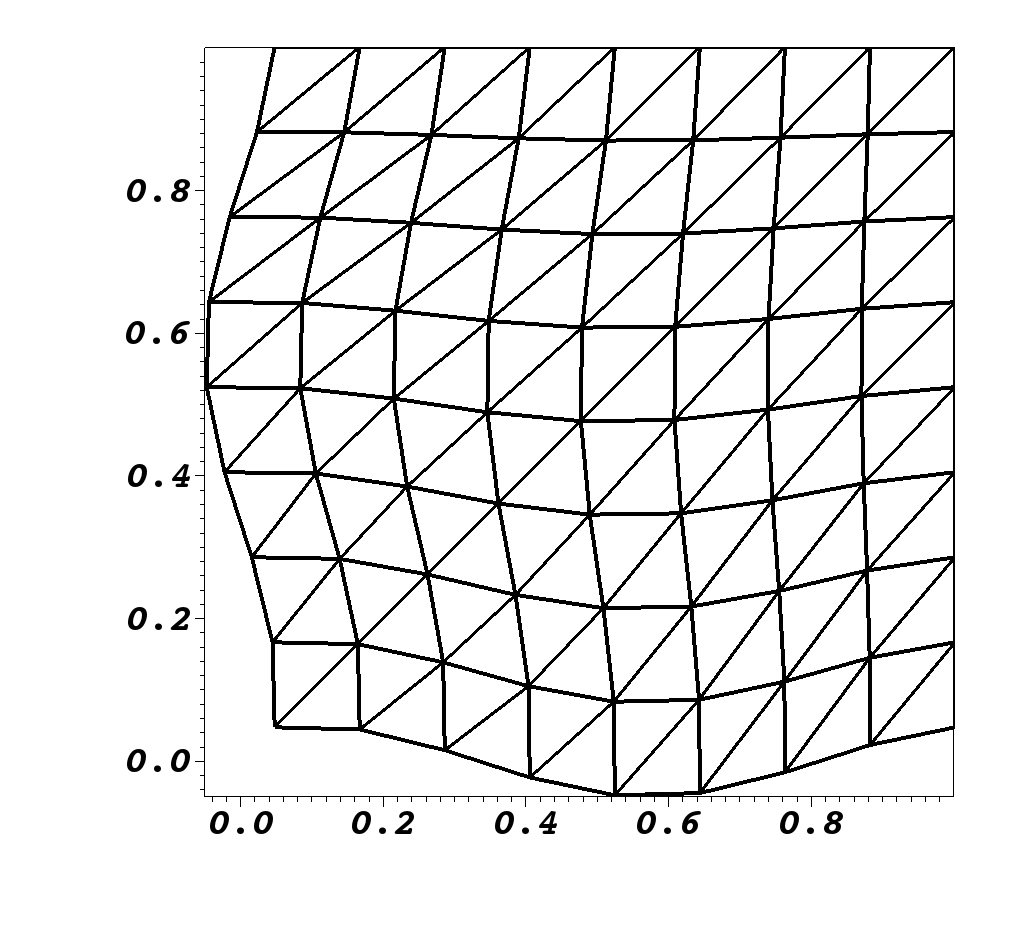}
    \caption{The mesh at different points in time for the test case
      of~\cref{ss:convrates}. From left to right the mesh at
      $t = 0, 0.4, 0.8$.\label{fig:mesh_movement}}
  \end{center}
\end{figure}

In this first test case we compute the rates of convergence of the
space--time HDG method applied to the Navier--Stokes equations on a
time-dependent domain. Introducing first a uniform triangular mesh for
the unit square, the mesh vertices $(x_1, x_2)$ for the deforming
domain $\Omega(t)$ are obtained at any time $t\in[0,1]$ by the
following relation
\begin{equation*}
  x_i = x_i^0 + 0.05 (1 - x_i^0)\sin( 2\pi (\tfrac{1}{2} - x_i^* + t) ) 
  \qquad  i = 1,2,  
\end{equation*}
where $(x_1^0, x_2^0) \in [0,1]^2$ are the vertices of the uniform
mesh and $(x_1^*, x_2^*) = (x_2^0, x_1^0)$. The mesh at three
different points in time is shown in~\cref{fig:mesh_movement}.

Let
$\partial\mathcal{E}^N := \cbr{(t,x_1,x_2)\in\partial\mathcal{E} :
  x_1=1}$
and
$\partial\mathcal{E}^D = \partial \mathcal{E} \setminus
(\partial\mathcal{E}^N \cup \Omega(0) \cup(\Omega(1))$.
The boundary conditions and source term $\boldsymbol{f}$ in~\cref{eq:navsto}
are chosen such that the exact solution is given by
\begin{equation*}
  \boldsymbol{u} = 
  \begin{bmatrix}
    (e^t - 1)\sin(\pi x_1) \sin(\pi x_2)\\
    (e^t - 1)\cos(\pi x_1) \cos(\pi x_2)
  \end{bmatrix},
  \qquad
  p = (2 + \cos(t))\sin(\pi x_1)\cos(\pi x_2).
\end{equation*}

We consider the rates of convergence for polynomial degrees $k=2$ and
$k=3$ and on a succession of refined space--time meshes. The coarsest
space--time mesh consists of $6 \cdot 8^2$ tetrahedra per space--time
slab with $\Delta t = 0.05$ and refinement happens in both space and
time. For the Picard iteration~\cref{eq:stopping_criterion} we set
$\text{TOL} = 10^{-12}$.

\Cref{tab:spatial_k2,tab:spatial_k3} show the rates of convergence for
the velocity and pressure computed on the spatial domain at final time
$T=1$ for $\nu = 10^{-4}$ and $\nu = 10^{-7}$, respectively. The rates
of convergence over the entire space--time domain $\mathcal{E}$ are
shown in~\cref{tab:spacetime_k2,tab:spacetime_k3} again for both
$\nu = 10^{-4}$ and $\nu = 10^{-7}$. We observe from these tables,
that for $k=2$ and for $k=3$, the method is optimal, i.e., the
velocity error is of order $\mathcal{O}(h^{k+1})$ and the pressure
error is of order $\mathcal{O}(h^k)$. The tables also show that the
error in the divergence of the approximate velocity is of machine
precision even on deforming domains, as expected
from~\cref{s:properties}.

Finally, we observe that the velocity error is independent of the
viscosity, i.e., the space--time hybridizable discontinuous Galerkin
method~\cref{eq:weakformulation_stns} is pressure-robust, even on
\emph{time-dependent} domains.

\begin{table}[tbp]
  \centering
  \begin{tabular}{ccccccc}
    \hline
    Cells per slab  & Nr. of slabs & $\norm{\boldsymbol{u} - \boldsymbol{u}_h}$  & rates  
    & $\norm{p - p_h}$ & rates & $\norm{\nabla \cdot \boldsymbol{u}_h}$ \\ 
    \hline
    \multicolumn{7}{l}{$\nu = 10^{-4}$} \\
    \hline
    $384$   & $20$  & 1.5e-3 &  -  & 1.8e-2 &  -  & 6.0e-14  \\ 
    $1536$  & $40$  & 2.0e-4 & 2.9 & 4.5e-3 & 2.0 & 1.1e-13  \\ 
    $6144$  & $80$  & 2.6e-5 & 2.9 & 1.1e-3 & 2.0 & 2.4e-13  \\ 
    $24576$ & $160$ & 3.3e-6 & 3.0 & 2.9e-4 & 1.9 & 4.7e-13  \\ 
    \hline
    \multicolumn{7}{l}{$\nu = 10^{-7}$} \\
    \hline
    $384$   & $20$  & 1.5e-3 &  -  & 1.8e-2 &  -  & 6.0e-14  \\ 
    $1536$  & $40$  & 2.0e-4 & 2.9 & 4.5e-3 & 2.0 & 1.2e-13  \\ 
    $6144$  & $80$  & 2.5e-5 & 3.0 & 1.2e-3 & 1.9 & 2.5e-13  \\ 
    $24576$ & $160$ & 3.3e-6 & 2.9 & 2.9e-4 & 2.0 & 5.0e-13  \\ 
    \hline    
  \end{tabular}
  \caption{The rates of convergence for the test case of~\cref{ss:convrates}
    with $k=2$ computed in the $L^2$-norm on $\Omega(T)$, the spatial domain 
    at final time $T=1$.}
  \label{tab:spatial_k2}
\end{table}

\begin{table}[tbp]
  \centering
  \begin{tabular}{ccccccc}
    \hline
    Cells per slab  & Nr. of slabs & $\norm{\boldsymbol{u} - \boldsymbol{u}_h}$  & rates  
    & $\norm{p - p_h}$ & rates & $\norm{\nabla \cdot \boldsymbol{u}_h}$ \\ 
    \hline
    \multicolumn{7}{l}{$\nu = 10^{-4}$} \\
    \hline
    $384$   & $20$  & 8.5e-5 &  -  & 1.3e-3 &  -  & 8.9e-13  \\ 
    $1536$  & $40$  & 5.2e-6 & 4.0 & 1.6e-4 & 3.0 & 1.8e-12  \\ 
    $6144$  & $80$  & 3.2e-7 & 4.0 & 2.0e-5 & 3.0 & 3.7e-12  \\ 
    $24576$ & $160$ & 2.0e-8 & 4.0 & 2.5e-6 & 3.0 & 7.6e-12  \\
    \hline
    \multicolumn{7}{l}{$\nu = 10^{-7}$} \\
    \hline
    $384$   & $20$  & 8.8e-5 &  -  & 1.3e-3 &  -  & 8.9e-13  \\ 
    $1536$  & $40$  & 5.6e-6 & 4.0 & 1.6e-4 & 3.0 & 1.9e-12  \\ 
    $6144$  & $80$  & 3.8e-7 & 3.9 & 2.0e-5 & 3.0 & 3.8e-12  \\ 
    $24576$ & $160$ & 2.8e-8 & 3.8 & 2.5e-6 & 3.0 & 7.6e-12  \\
    \hline    
  \end{tabular}
  \caption{The rates of convergence for the test case of~\cref{ss:convrates}
    with $k=3$ computed in the $L^2$-norm on $\Omega(T)$, the spatial domain 
    at final time $T=1$.}
  \label{tab:spatial_k3}
\end{table}

\begin{table}[tbp]
  \centering
  \begin{tabular}{ccccccc}
    \hline
    Cells per slab  & Nr. of slabs & $\norm{\boldsymbol{u} - \boldsymbol{u}_h}$  & rates  
    & $\norm{p - p_h}$ & rates & $\norm{\nabla \cdot \boldsymbol{u}_h}$ \\ 
    \hline
    \multicolumn{7}{l}{$\nu = 10^{-4}$} \\
    \hline
    $384$   & $20$  & 6.9e-4 &  -  & 6.8e-3 &  -  & 4.2e-14 \\ 
    $1536$  & $40$  & 8.8e-5 & 3.0 & 1.7e-3 & 2.0 & 8.3e-14 \\ 
    $6144$  & $80$  & 1.2e-5 & 2.9 & 4.3e-4 & 2.0 & 1.6e-13 \\ 
    $24576$ & $160$ & 1.5e-6 & 3.0 & 1.1e-4 & 2.0 & 3.0e-13 \\ 
    \hline
    \multicolumn{7}{l}{$\nu = 10^{-7}$} \\
    \hline
    $384$   & $20$  & 6.9e-4 &  -  & 6.8e-3 &  -  & 4.3e-14 \\ 
    $1536$  & $40$  & 8.6e-5 & 2.9 & 1.7e-3 & 2.0 & 8.9e-14 \\ 
    $6144$  & $80$  & 1.1e-5 & 2.9 & 4.3e-4 & 2.0 & 1.8e-13 \\ 
    $24576$ & $160$ & 1.4e-6 & 2.9 & 1.1e-4 & 2.0 & 3.5e-13 \\ 
    \hline    
  \end{tabular}
  \caption{The rates of convergence for the test case of~\cref{ss:convrates}
    with $k=2$ computed in the $L^2$-norm on $\mathcal{E}$, the whole space--time domain.}
  \label{tab:spacetime_k2}
\end{table}

\begin{table}[tbp]
  \centering
  \begin{tabular}{ccccccc}
    \hline
    Cells per slab  & Nr. of slabs & $\norm{\boldsymbol{u} - \boldsymbol{u}_h}$  & rates  
    & $\norm{p - p_h}$ & rates & $\norm{\nabla \cdot \boldsymbol{u}_h}$ \\ 
    \hline
    \multicolumn{7}{l}{$\nu = 10^{-4}$} \\
    \hline
    $384$   & $20$  & 3.4e-5 &  -  & 3.4e-4 &  -  & 6.8e-13 \\ 
    $1536$  & $40$  & 2.1e-6 & 4.0 & 4.3e-5 & 3.0 & 1.4e-12 \\ 
    $6144$  & $80$  & 1.2e-7 & 4.1 & 5.4e-6 & 3.0 & 2.7e-12 \\ 
    $24576$ & $160$ & 7.5e-9 & 4.0 & 6.8e-7 & 3.0 & 5.5e-12  \\
    \hline
    \multicolumn{7}{l}{$\nu = 10^{-7}$} \\
    \hline
    $384$   & $20$  & 3.5e-5 &  -  & 3.4e-4 &  -  & 6.8e-13 \\ 
    $1536$  & $40$  & 2.3e-6 & 3.9 & 4.3e-5 & 3.0 & 1.4e-12 \\ 
    $6144$  & $80$  & 1.5e-7 & 3.9 & 5.4e-6 & 3.0 & 2.8e-12 \\ 
    $24576$ & $160$ & 1.1e-8 & 3.8 & 6.8e-7 & 3.0 & 5.6e-12  \\
    \hline    
  \end{tabular}
  \caption{The rates of convergence for the test case of~\cref{ss:convrates}
    with $k=3$ computed in the $L^2$-norm on $\mathcal{E}$, the whole space--time domain.}
  \label{tab:spacetime_k3}
\end{table}

\subsection{Flow around a cylinder}
\label{ss:flowcylinder}

Next we consider flow around a cylinder. The setup of this test case
is taken from~\cite{Schaefer:1996, Lehrenfeld:2016} in which we
consider a fixed spatial domain $[0; 2.2]\times [0; 0.41]$ with a
cylindrical obstacle with radius $r = 0.05$ centred at
$(x_1,x_2) = (0.2, 0.2)$. 

\begin{figure}[tbp]
  \begin{center}
    \includegraphics[width=\linewidth]{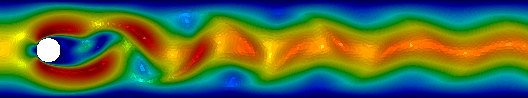}
    \caption{The velocity magnitude of flow around a cylinder, as
      described in~\cref{ss:flowcylinder}, at $t = 5$ using $3222$
      triangles and $k = 3$.\label{fig:cyl_final}}
  \end{center}
\end{figure}

A homogeneous Neumann boundary condition is applied on the outflow
boundary at $x_1=2.2$. On the inflow boundary at $x_1=0$ we impose
$\boldsymbol{u} = [6x_2(0.41-x_2)/0.41^2, 0]^T$, while $\boldsymbol{u} = \boldsymbol{0}$ is
imposed on the cylinder and on the walls $x_2=0$ and $x_2=0.41$. The
kinematic viscosity is set to be $\nu = 10^{-3}$ and the initial
condition is obtained by solving the steady--state Stokes problem, as
done also in~\cite{Rhebergen:2018a}. For the stopping criterion of the
Picard iteration we used $\text{TOL} = 10^{-9}$. The velocity
magnitude at final time $t=5$ is shown in~\cref{fig:cyl_final}.

To validate the dicretization we compute the lift and drag
coefficients. These are defined as
\begin{equation*}
C_L = \frac{1}{r \Delta t} \int_{\Gamma_c} (\sigma_d \cdot \boldsymbol{n})\cdot \boldsymbol{e}_1, \qquad 
C_D = \frac{1}{r \Delta t} \int_{\Gamma_c} (\sigma_d \cdot \boldsymbol{n})\cdot \boldsymbol{e}_2,
\end{equation*}
where $\boldsymbol{e}_1$ and $\boldsymbol{e}_2$ are the unit vectors in the $x_1$
and $x_2$ directions, respectively, and $\Gamma_c$ denotes the
space--time boundary of the cylinder.

We compute the lift and drag coefficients on a space--time mesh
consisting of $9666$ tetrahedra per slab. We consider a time step of
$\Delta t = 5\cdot 10^{-3}$. Over the time interval $t\in [0,5]$ we
found that $C_L \in [-1.014, 0.98]$ and $C_D \in [3.153, 3.219]$,
which compare well to the results found in
literature~\cite{Schaefer:1996, Lehrenfeld:2016}.

\subsection{Flow around an oscillating airfoil}
\label{ss:airfoil}

In this final test case we consider the simulation of flow around an
oscillating NACA0012 airfoil on the domain $[-5,10] \times [-5, 5]$
with trailing edge at the origin. The computational domain consists of
17088 tetrahedra per space--time slab and as polynomial approximation
we use $k=2$. We set the kinematic viscosity to be $\nu = 10^{-3}$. As
stopping criterion in the Picard iteration we use
$\text{TOL} = 10^{-6}$.

\begin{figure}[tbp]
  \begin{center}
    \subfloat[Mesh at $10^{\circ}$ angle of attack.]{\includegraphics[width=.45\textwidth]{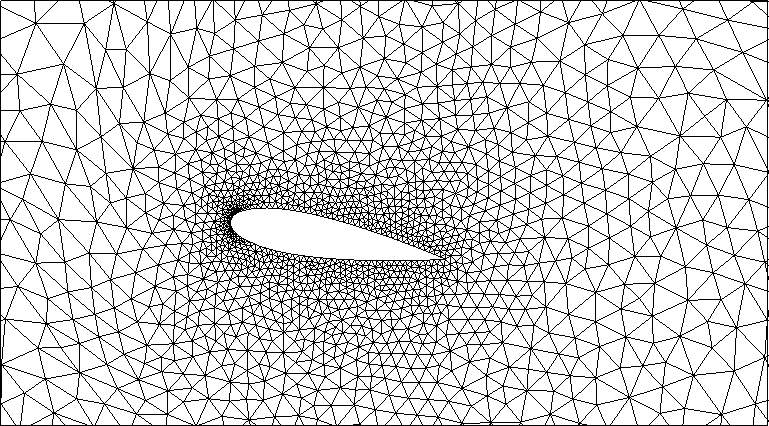}} \quad
    \subfloat[Mesh at $30^{\circ}$ angle of attack.]{\includegraphics[width=.45\textwidth]{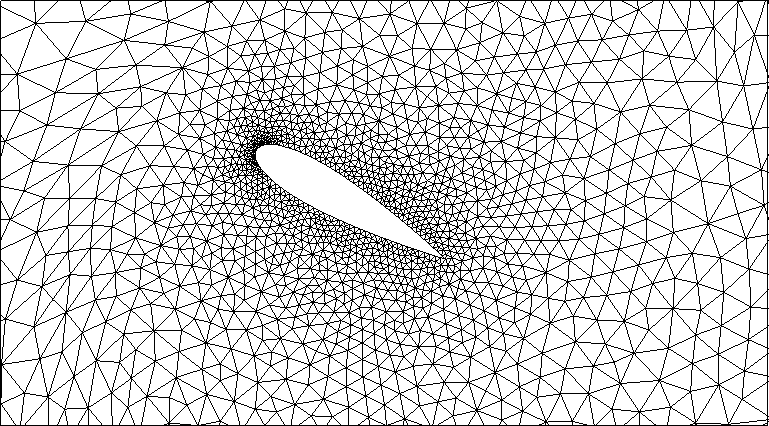}}
    \caption{Mesh at different angles of attack for the test case
      described in~\cref{ss:airfoil}. \label{fig:naca_meshes}}
  \end{center}
\end{figure}

To obtain the initial condition we first solve the steady Stokes
problem around the airfoil at an angle of attack of $20^{\circ}$. We
then solve the Navier--Stokes problem around the fixed airfoil from
$t=0$ to $t=1$ using a time step of $\Delta t = 0.01$. Given this
`initial condition' we then prescribe an oscillatory movement of the
airfoil for $t>1$. Keeping the trailing edge fixed at
$(x_1,x_2)=(0,0)$, the angle of attack changes according to
\begin{equation*}
  \text{angle of attack} = 20 + 10\sin(\pi (t-1)).
\end{equation*}
For $t>1$ we use $\Delta t = 2\cdot 10^{-3}$. To account for the time
dependent angle of attack, the mesh is updated at each time step as
follows. Nodes within a radius of $1.5$ from the trailing edge move
with the airfoil, nodes outside a radius of $2$ from the trailing edge
remain fixed, while the movement of the remaining nodes decrease
linearly with distance, see~\cref{fig:naca_meshes}.

\begin{figure}[h]
  \begin{center}
    \subfloat[Airfoil at $20^{\circ}$ angle of attack at $t=1$.]{\includegraphics[width=.45\textwidth]{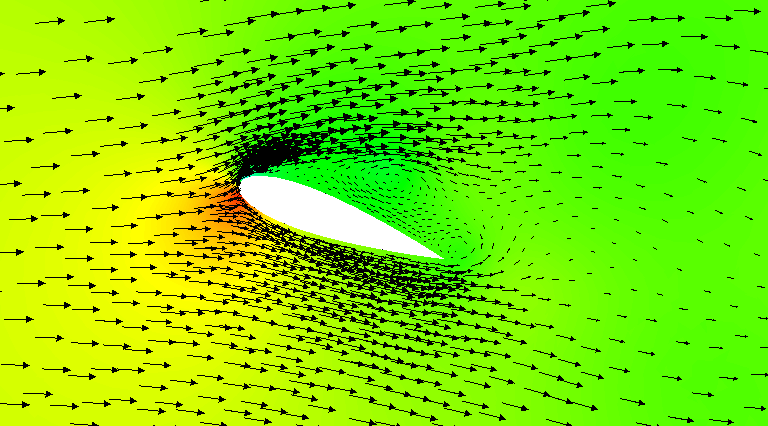}} \quad
    \subfloat[Airfoil at $30^{\circ}$ angle of attack at $t=1.5$.]{\includegraphics[width=.45\textwidth]{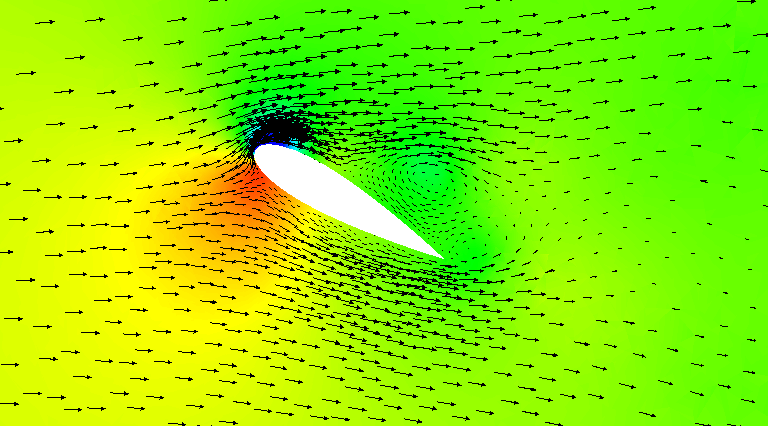}} \\
    \subfloat[Airfoil at $20^{\circ}$ angle of attack at $t=2$.]{\includegraphics[width=.45\textwidth]{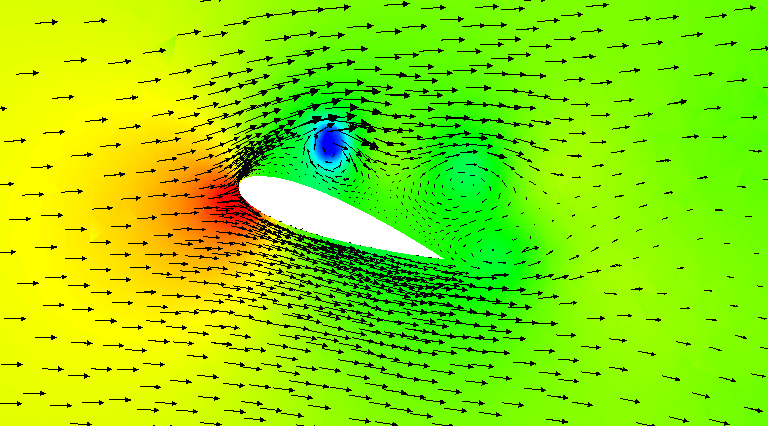}} \quad
    \subfloat[Airfoil at $10^{\circ}$ angle of attack at $t=2.5$.]{\includegraphics[width=.45\textwidth]{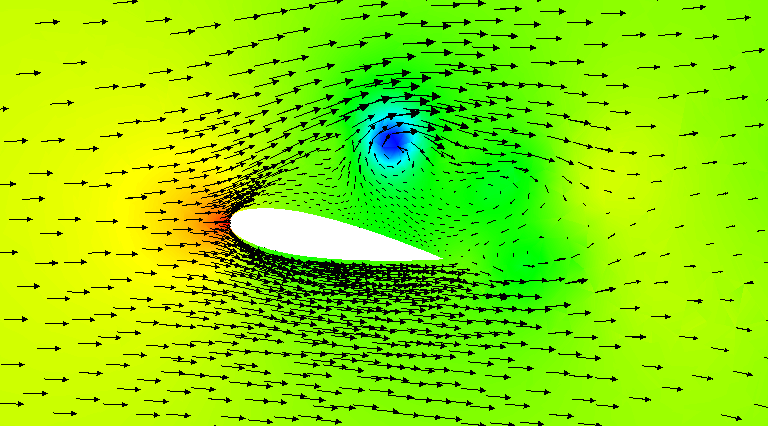}} \\
    \subfloat[Airfoil at $20^{\circ}$ angle of attack at $t=3$.]{\includegraphics[width=.45\textwidth]{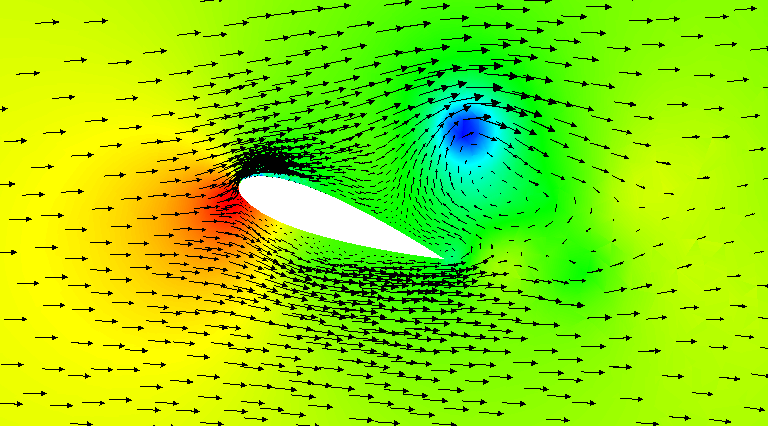}} \quad
    \subfloat[Airfoil at $30^{\circ}$ angle of attack at $t=3.5$.]{\includegraphics[width=.45\textwidth]{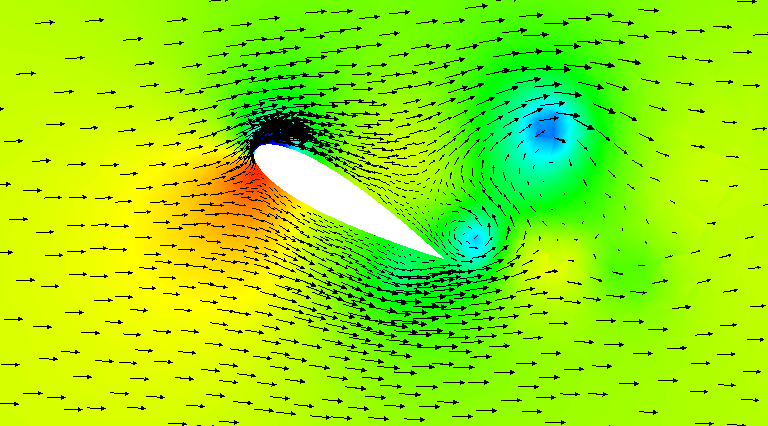}} \\
    \subfloat[Airfoil at $20^{\circ}$ angle of attack at $t=4$.]{\includegraphics[width=.45\textwidth]{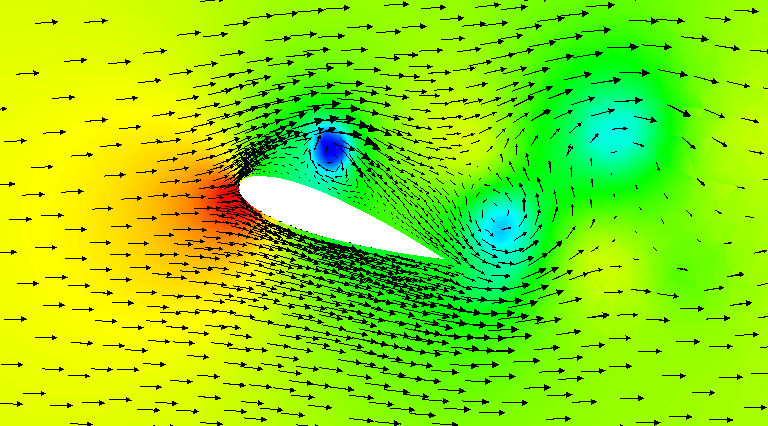}}
    \caption{Velocity vector plot and pressure field around an
      oscillating NACA0012 airfoil as described
      in~\cref{ss:airfoil}.\label{fig:naca}}
  \end{center}
\end{figure}

In~\cref{fig:naca} we plot the computed pressure and velocity vector
fields. We observe vortex shedding at the trailing edge and detachment
of vortices over the top of the airfoil while new vortices form at the
tip of the airfoil. These phenomena agree with those observed in
literature, for example~\cite{Tezduyar:1992}.

\section{Conclusions}
\label{s:conclusions}

We presented a space--time hybridizable discontinuous Galerkin finite
element method for the Navier--Stokes equations on time-dependent
domains. This scheme guarantees a point-wise divergence-free and
$H({\rm div})$-conforming velocity field, it is locally momentum
conserving and energy stable, even on dynamic meshes. We have shown
the performance of the method in terms of rates of convergence,
pressure-robustness, and flow simulations around a cylinder and an
oscillating airfoil.

\section*{Acknowledgments}
SR gratefully acknowledges support from the Natural Sciences and
Engineering Research Council of Canada through the Discovery Grant
program (RGPIN-05606-2015) and the Discovery Accelerator Supplement
(RGPAS-478018-2015).

\bibliographystyle{abbrvnat}
\bibliography{references}
\end{document}